\documentclass[12pt]{article}
\usepackage{amsthm, amsfonts, url}
\usepackage{amscd,amsmath,amssymb}

\newtheorem{theorem}{\noindent Theorem}
\newtheorem{lemma}{\noindent Lemma}
\newtheorem{definition}{\noindent Definition}
\newtheorem{corollary}{\noindent Corollary}

\newtheorem{proposition}{\noindent Proposition}

\urldef\vershik\url{vershik@pdmi.ras.ru}
\title{Towards the definition of metric hyperbolicity}

\author{A.~M.~Vershik\thanks{%
St.~Petersburg Department of Steklov Institute of Mathematics.
E-mail: \vershik.
 Partially supported by grants NSh-2251.2003.1, CRDF RUM1-2262-ST-04,
and INTAS 03-51-5018.}}

\begin{document}

\maketitle

\rightline{\it To Ya.~G.~Sinai for his 70th anniversary.}

\smallskip
\begin{abstract}
We introduce measure-theoretic definitions of {\it hyperbolic
structure for mea\-su\-re-preserving automorphisms}. A wide class of
$K$-auto\-morph\-isms possesses a hyperbolic structure; we prove that
all $K$-auto\-morph\-isms have a slightly weaker structure of {\it
semi-hyperbolicity}. Instead of the notions of stable and unstable
foliations and other notions from smooth theory, we use the tools
of the theory of polymorphisms. The central role is played by  {\it
polymorphisms} associated with a special invariant equivalence relation,
more exactly,  with a
homoclinic equivalence relation. We call an
automorphism with given hyperbolic structure a hyperbolic
automorphism
and prove that it is canonically
quasi-similar to a so-called prime nonmixing polymorphism.
 We present a short
but necessary vocabulary of polymorphisms and Markov operators
from \cite{V1,V2}.
\end{abstract}

\tableofcontents

\section{Motivations and statement of the problem}

The theory of hyperbolic dynamical systems is one of the main
achievements in the theory of dynamical systems of the second half of the
last century. Although the basic concept appeared as far ago as
in the papers by H.~Poincar\'e and J.~Hadamard and was discussed in many
subsequent papers, and the main example --- the geodesic flow on a
surface of constant negative curvature --- was known from the
very beginning, and some ``hyperbolic'' effects (such as exponential rate of
convergence and divergence, Lyapunov exponents, etc.)
were known in terms of
concrete differential equations (such as the Van der Pol equation, which
was studied by G.~Littlewood and M.~Cartwright), but in the framework of the
modern theory of dynamical systems, ergodic theory, and
representation theory, it was considered only  in the 40s--50s by E.~Hopf and
G.~Hedlund and I.~M.~Gelfand and S.~V.~Fomin. The analysis of
concrete examples gave an impulse to the general theory, which was
formulated and axiomatized
in the 60s by several authors (S.~Smale,
D.~Anosov, and others). The definition of smooth hyperbolic systems
involves the notions of Riemannian metric, stable and unstable
foliations on a manifold, etc., which use the smooth structure of
the phase space.

At the same time, connections of these ideas with the theory of
stationary random processes were advocated already in the 40s by
A.~N.~Kolmogorov, who considered this type of dynamical systems
in a very wide context; he defined the notion of regular
random stationary processes; apparently, he was the first to
emphasize that the sigma-fields of the partitions
with fixed ``past''
and ``future'' of stationary random processes are similar to
pairs of horocycle foliations or geodesic flows on the unit
tangent bundle of compact surfaces of constant negative
curvature. In 1958, Kolmogorov introduced a Shannon-type entropy
{\it as a metric invariant of dynamical systems} and solved the
long-standing problem on isomorphisms of Bernoulli systems. His
pupil Ya.~Sinai, together with V.~A.~Rokhlin and their schools,
developed, in the 60s--70s, entropy theory and the theory of
$K$-systems in the framework of ergodic theory and dynamical systems.
Ya.~Sinai's contributions concerned not only the theory of dynamical
systems, but much wider areas, including statistical
physics, classical dynamics, ergodic hypothesis, and so on; his
activity helped to combine dynamical theory with statistical
physics and many others topics.

 The link between classical hyperbolic systems and the class of
 $K$-systems became more clear in the 70s, after the
 papers by D.~Ornstein appeared,
 who, starting from Sinai's theorem on weak isomorphisms of
 Bernoulli automorphisms with the same entropy, proved a fundamental
 result that entropy is a complete invariant in the class of Bernoulli
 systems and gave an invariant definition of Bernoulli systems.
 This result allowed him together with B.~Weiss \cite{OW} to prove the
 Bernoulli property of the geodesic flow on a compact surface of
 constant negative curvature. Bernoulli property of the hyperbolic
 automorphisms of torus also follows from that theory.
 Later, Ya.~Pesin \cite{Pes} proved the
 Bernoulli property  for smooth hyperbolic systems in full generality.
 The existence of non-Bernoulli $K$-systems, which was discovered by
 Ornstein and Shields, and especially Kalikow's example opened a
 new class of problems in dynamical theory.

 But we can see a gap between hyperbolicity in the smooth
 category and $K$-property in the category of measure spaces --- we have no
 purely measure-theoretic analogs of the notions of hyperbolic
 theory. The vague analogy
 between stable and unstable foliations on one hand
 and the ``past'' and ``future'' of a stationary $K$-process
 on the other have not been put
 into an appropriate general scheme.

 To be more concrete, let us formulate our main problem:

\smallskip
{\bf The goal of this paper is to suggest a purely measure-theoretic definition of
  hyperbolic structure of measure-preserving automorphisms and to develop
  some tools for studying it.}
\smallskip

  In order to do this, we must overcome difficulties with
  definitions of objects that use the smooth and metric structures
  and the corresponding numerical characteristics.

  To this end, we use the notion of polymorphism
  (= Markov, or multivalued, map in ergodic theory, see \cite{V1})
  and the corresponding tools.
  This allows us to avoid problems with defining
  foliations and so on: roughly speaking, instead of foliations
  we consider polymorphisms which are associated with them. This
  allows us to transfer metrical and topological questions, including
  estimations, from the manifold (phase space) to the space of transformations or
  operators. In order to formulate the main definition,
  we use the weak topology in the space of polymorphisms
  and Markov operators
  (Condition \textbf{H} and the definition of hyperbolicity, see Sec.~3).

   Roughly speaking, a hyperbolic structure for an automorphism with invariant measure
   is a polymorphism associated with an invariant equivalence relation
   which plays the role of the homoclinic equivalence. We called such automorphisms
    hyperbolic automorphisms; the definition is metrically invariant.
   An automorphism
   can have several hyperbolic structures or none at all.
   Presumably, a hyperbolic automorphism must be a
   $K$-automorphism satisfying an additional property
   (property ($*$), see Sec.~3.3 and the problem in
   Sec.~3.7). We also define
    a weaker notion of semi-hyperbolic structure, which could be related to
    the notion
    of partial hyperbolic systems in the sense of Pesin (see \cite{Pes}).

    An extremely important notion closely related to this topic is the notion
    of quasi-similarity; one of our main results claims that the
    hyperbolicity of an automorphism $T$ implies its quasi-similarity with a
    prime nonmixing and non-co-mixing polymorphism (see Sec.~3.2). The notion of
    quasi-similarity came from the theory of contractions in Hilbert spaces
    (see \cite{FN}) and scattering theory (\cite{Lax}) and was not used
    earlier in the theory of dynamical systems. One of the motivations
    of this paper is to study the interrelation between classical dynamical systems
    and polymorphisms and apply it to hyperbolic theory; in particular,
    to investigate the notion of quasi-similarity between polymorphisms and automorphisms.
    The theory of polymorphisms and Markov operators has also a direct contact
    with the theory of Markov processes, which can be used for refining some of our
    results. We discussed this question briefly in \cite{V1,V2}.

   We start in Sec.~2 with briefly recalling some of the notions
   concerning polymorphisms
   and Markov operators. Section~3 contains out main results. We give
   definitions and first corollaries of hyperbolicity in Sec.~3.1
   and a theorem on quasi-similarity in Sec.~3.2; explain
   how to include classical examples
   into our approach in Sec.~3.3; present a
   geometrical interpretation of quasi-similarity in Sec.~3.4; give
   the definition of semi-hyperbolic structures in Sec.~3.5;
   and prove that all $K$-automorphisms
   have a semi-hyperbolic structure in Sec.~3.6. Two questions from
   a large list of open problems are presented in Sec.~3.7.

The author is grateful
to N.~Tsilevich for her help with the preparation of the manuscript.

\section{Vocabulary of polymorphisms and Markov operators}
\subsection{Polymorphisms}
    In order to make this paper independent, we give a short list
    of definitions concerning polymorphisms and Markov operators.
    The reader can find a detailed version in \cite{V1}. In other
    areas of mathematics, notions parallel to that of polymorphism
    are: {\it correspondence} in algebra and algebraic geometry;
    {\it bifibration} in differential geometry, {\it Markov map} in probability
    theory, {\it Young measure} in optimal control, etc.
    Equivalent definitions of the notions under consideration
    in terms of Markov operators will be presented in the next section.

 \begin{definition}
    {A polymorphism with invariant measure $\Pi$
    of a Lebesgue space $(X,\mu)$ to itself is a diagram consisting
    of an ordered triple of Lebesgue spaces:
    $$ (X,\mu)  \stackrel{\pi_1}\longleftarrow (X \times X, \nu)
    \stackrel{\pi_2}\longrightarrow (X, \mu),$$
    where $\pi_1$ and $\pi_2$ stand for the projections to the first
    and second component of the product space $(X \times X, \nu)$, and
    the measure $\nu$, which is defined on the $\sigma$-field generated by the
    product of the $\sigma$-fields of $\bmod0$ classes of measurable
    sets in $X$, is such that $\pi_i\nu=\mu$, $i=1,2$.
    The measure $\nu$ is called the {\it bistochastic measure of the polymorphism $\Pi$}.}
\end{definition}

     Let us define the main structures on the set of the polymorphisms and
     the notions of the theory of the polymorphisms.

    1.A polymorphism $\Pi^*$ is called {\it conjugate} to the polymorphism
    $\Pi$ if its diagram is obtained from the diagram of  $\Pi$
    by reflecting with respect to the central term. If
    polymorphism is an automorphism then conjugate polymorphism is nothing more than
    inverse automorphism.

    Consider the ``vertical'' partition $\xi_1$ and
    the ``horizontal'' partition $\xi_2$ of the space $(X\times X, \nu)$
    into the preimages of points under the projections $\pi_1$ and $\pi_2$,
    respectively.
    In terms of bistochastic measures, the
    {\it value of a polymorphism at a point $x\in X$}
    is a conditional measure. More precisely, we have the following definition.

\begin{definition}
    {In the above notation, the {\it value} of the polymorphism
    $\Pi:X \to X $ at a point $x_1 \in X$
    is, by definition,  the conditional measure $\nu^{x_1}$
    of $\nu$ on the set $\{(x_1,\cdot)\}$ with respect to the vertical
    partition $\xi_1$ (the transition probability); similarly,
    the {\it value} of the conjugate polymorphism $\Pi^*$ at a point $x_2\in X$
    is the conditional measure $\nu_{x_2}$ of $\nu$ on the set $\{(\cdot,x_2)\}$
    with respect to the horizontal partition $\xi_2$ (the cotransition probability).
    These systems of conditional measures are well-defined on sets of full $\nu$-measure.
    Thus a polymorphism is a ${\mod 0}$ class of measurable maps from $(X,\mu)$ to the
    space of measures on $X$ equipped with the ordinary Borel structure.}
\end{definition}

\smallskip

    2. The set of {\it ``images''} of the points $x\in X$
    under a polymorphism is the system of conditional (transition)
    measures $\nu^x$ on $X$. It is
    important that this system is defined up to measure zero, so
    there is no sense in the measure-theoretic category of ``individual'' image,
    but only the system of images as a whole makes sense.
    If almost all measures $\nu^x$ are delta-measures, then
    we have a deterministic measure-preserving map.
    Denote the set of polymorphisms of a given space $(X,\mu)$
    by $\mathcal P(X,\mu)=\mathcal P$.

\smallskip

     3. {\it A multiplication} in the set of polymorphisms (bistochastic
    measures) $\mathcal P$ is defined as follows:
    let $\Pi_1,\Pi_2$ be two polymorphisms
    with bistochastic measures $\nu_1,\nu_2$; then the
    product $\Pi_1\Pi_2$  has the bistochastic measure $\nu$ defined by
    $$\nu^x(A)=\int \nu_1^y(A)\,  d\nu_2^x(y).$$
    The ordinary weak topology on the set of polymorphisms (a neighborhood of the identity
    is the set of polymorphisms whose bistochastic measures are in a neighborhood
    of the diagonal measure) endows it with the structure
    of a {\it compact topological semigroup}. This semigroup has a unity (the identity map),
    involution (conjugacy), zero element $\Theta$ (the bistochastic product
    measure $\nu=\mu \times \mu $), and a {\it natural convex structure}
    on the set of bistochastic measures.
    The subgroup of invertible elements of the semigroup $\cal P$
    is the group of measure-preserving transformations.
    The set of all polymorphisms of a finite space is the
    convex semigroup of bistochastic matrices.
    The convex semigroup of all polymorphisms of a Lebesgue space with
    a continuous measure is the inverse limit of the sequence of convex
    compact spaces of bistochastic matrices.

\smallskip
     Now define the classes of polymorphisms, factor polymorphisms, ergodicity, mixing,
     primeness and quasisimilarity of polymorphisms.

     4. A measurable partition $\xi$ is called {\it invariant} under a
     polymorphism $\Pi$ if for almost all elements $C \in \xi$ there exists
     another element $D \in \xi$ such that for almost all (with respect to the
     conditional measure on $C$) points $x \in C$, we have $\mu^x(D)=1$, where $\mu^x$
     is the $\Pi$-image of $x$.
     In other words, the factor polymorphism $\Pi_{\xi}$ of $\Pi$
     with respect to an invariant partition $\xi$ is an endomorphism of
     the space $(X_{\xi},\mu_{\xi})$.

     In particular, if for almost all elements $C \in \xi$ of
     $\xi$ we have $\mu^x(C)=1$ for almost all $x\in C$, then
     $\xi$ is called a {\it fixed} partition for $\Pi$
     and the corresponding factor polymorphism is the identity map on $X_\xi$.
     A polymorphism is {\it ergodic} if it has no nontrivial identity
     factor automorphism.

\smallskip

    5. For a given polymorphism $\Pi$, of
    a space $(X,\mu)$, with bistochastic measure $\nu$,
    the {\it factor polymorphism} (or {\it quotient}) of $\Pi$
    by a measurable partition $\xi$
    is the polymorphism of the space $(X/{\xi},\mu/{\xi})$ to itself
    with the factorized bistochastic measure $\nu/{(\xi \times \xi)}$.
    Thus the factor polymorphism of any polymorphism by any measurable
    partition does exist; in particular, the {\it factor polymorphism of
    any automorphism by any (not necessarily invariant) partition
    always does exist}.

    A polymorphism is called {\it prime} if it has no nontrivial
    invariant partition, or has no nonzero factor endomorphism.
    A polymorphism is called {\it coprime} if its conjugate is prime.
    (Compare this notions with the notions of a prime
    automorphism and an exact endomorphism.)

\smallskip

  6. A polymorphism $\Pi$ is called {\it mixing} if the sequence of its powers
     tends to the zero polymorphism $\Theta$ in the weak topology:
     \mbox{w-$\lim_{n \to \infty}\Pi^n = \Theta$}.  Note that
     it may happen that a polymorphism is mixing while its conjugate is not.
     We call a polymorphism {\it co-mixing} if its conjugate is a
     mixing polymorphism.

\smallskip

  7. A polymorphism $\Pi$ is called {\it semi-dense} if a measurable function
     for which $\int f(y)d\nu^x(y)=0$ for $\mu$-almost all $x$ is equal
     to zero; a polymorphism $\Pi$ is {\it dense} if
     both $\Pi$ and
     $\Pi^*$ are semi-dense. It is more convenient to express density in terms
     of Markov operators (see the next subsection).

\smallskip

  8. A polymorphism $\Pi$ is called {\it nondegenerate}
  if for almost all $x$, the conditional measure $\nu^x$ of $\Pi$
  is not a delta-measure.

\smallskip\noindent
\textbf{Remark}. Sometimes it is more convenient to regard a
partition as an equivalence relation (congruence); we will not
distinguish an equivalence relation on $(X, \mu)$ and
the corresponding partition of $(X, \mu)$, and will denote them by the same
letter.

\smallskip

9. We say that a polymorphism $\Pi$ {\it is associated with a partition}
$\xi$ (measurable
  or not) if for almost all $x$ we have $\nu^x(\xi(x))=1,$ where $\xi(x)$ is
the element
  of $\xi$ that contains $x$. In other words, the polymorphism
  acts {\it along the blocks of the partition}. Each automorphism
  is associated with its orbit partition (see \cite{V1}).

\smallskip

  10. A polymorphism (in particular, automorphism) $\Pi$ is a {\it
  quasi-image} (see an analog of this notion in \cite{FN})
  of a polymorphism or automorphism $\Psi$
  if there exists a {\it dense} polymorphism $\Lambda$ such that
  $\Lambda \Pi =\Psi \Lambda$. If $\Psi$ is also
  a quasi-image of $\Pi$, then we say that $\Pi$ and $\Psi$ are {\it quasi-similar}.
  Quasi-similarity is a much more rough equivalence than
  similarity, for example,
  mixing is not an invariant of quasi-similarity.
  The classification of automorphisms (e.g., of $K$-automorphisms) up to
  quasi-similarity is a very intriguing problem; one problem of such a type:
  is entropy of automorphism an invariant under quasisimilarity?
   But for further discussions it is especially important that an automorphism may be
  quasi-similar to a polymorphism; our main definition uses this
  fact in a very crucial way (see the definition of hyperbolicity).

    \subsection{Markov operators}

   The functional analog of the notion of  polymorphism is the notion of
   {\it Markov operator} in the Hilbert space $L^2_{\mu}(X)$ (see the
   classical theory in \cite{R, V2}).

\begin{definition}
{\it A Markov operator in the Hilbert space $L^2(X,\mu)$ of
   complex-valued square integrable functions on a Lebesgue--Rokhlin
   space $(X,\mu)$ with a continuous normalized measure $\mu$ is a
   continuous linear operator $V$ satisfying the following
   conditions:

{\rm 1)} $V$ is a contraction: $\|V\|\leq 1$ (in the operator
norm);

{\rm 2)} $V{\emph{1}}= {V}^*{\emph{1}}={\emph{1}}$, where
$\emph{1}$
   is the function identically equal to one;

{\rm 3)} $V$ preserves the nonnegativity of functions: $Vf$ is
   nonnegative whenever $f\in L^2(X,\mu)$ is nonnegative}.
\end{definition}

   Note that condition 1) follows from  2) and 3), and the second
   condition in 2) follows from the other ones.
   In short: a Markov operator is a unity-preserving positive
   contraction.

  The set $\mathcal M$
  of all Markov operators is a convex weakly compact semigroup with
  involution $V \to V^*$. Unitary (isometric) Markov operators are precisely
  the operators generated by measure-preserving
  auto(endo)morphisms.

\begin{proposition}
  {\rm 1.} Let $\Pi$ be a polymorphism of a space $(X, \mu)$ with
  invariant measure; then the formula
  $$(V_{\Pi} f)(x)= \int_X f(y)\mu^x(dy)$$
  correctly defines a Markov operator in $L^2_\mu(X)$.

  {\rm 2.} Every Markov operator $V$ in the space $L^2_\mu(X)$, where
  $(X,\mu)$ is a Lebesgue space with continuous finite measure,
  can be represented in the form $V=V_{\Pi}$, where $\Pi$ is
  a  polymorphism of $(X,\mu)$ with invariant measure.

  {\rm 3.} The correspondence $\Pi \mapsto V_{\Pi}$ is a
  continuous (with respect to the weak topologies) antiisomorphism between the
  convex compact semigroup with involution of $\bmod 0$ classes of
  polymorphisms and the analogous semigroup of Markov operators.
\end{proposition}

 The Markov operator $\textbf{1}=W_{\Theta}$ corresponding to the zero polymorphism
 $\Theta$ is the orthogonal projection to the one-dimensional subspace
 of constants. The operator of mathematical expectation is, obviously,
 also a Markov operator;  it corresponds to the polymorphism
 that sends a point $x$ to the conditional measure of the element
 of the partition (corresponding to the expectation) that contains $x$.

\smallskip
 Now we reformulate the notions introduced for polymorphisms (ergodicity,
 mixing, primality, density, etc.) in terms of Markov operators.

\smallskip
 1. A Markov operator $V$ is called {\it mixing}
 (resp. {\it comixing}) if the sequence $V^n$
 (resp. ${V^*}^n$) weakly tends, as $n \to \infty$, to the projection onto the
 subspace of constants: $$V^n \to \textbf{1}=V_{\Theta},
 (\mbox{resp.}{V^*}^n \to \textbf{1}=V_{\Theta}).$$
 The Markov operator $V_{\Pi}$ is mixing (comixing) if and only if the polymorphism
 $\Pi$ is mixing (resp. comixing).

\smallskip
 2. We will say that a Markov operator $V=V_{\Pi}$ is {\it semi-dense} if the
 $V$-image of the space $L^2_{\mu}(X)$ is dense in $L^2_{\mu}(X)$;
 this is equivalent to the triviality of the kernel of the conjugate
 operator and, consequently, to the semi-density of the polymorphism $\Pi$. A
 Markov operator $V$ is called {\it dense} if both $\ker V$ and
 $\ker V^*$ are trivial. The density of $V_{\Pi}$ is equivalent to
 the density of $\Pi$.

 \smallskip
 3. A Markov operator $V$ is a {\it quasi-image} of a Markov operator $W$ if
 there exists a semi-dense Markov operator $L$ such that $ WL=LV$.
 Two Markov operators are {\it quasi-similar} if each of them is a
 quasi-image of the other one. Two Markov operators are quasisimilar
 if and only if the corresponding polymorphisms are quasisimilar.

 \smallskip
 4. A Markov operator $V$ is called {\it totally nonisometric} if
 there is no nonzero invariant subalgebra\footnote{More exactly, a subspace that
 consists of all functions from $L^2$ that are constant a.e.\
 on all elements of some measurable partition, see \cite{Rokh}.}
 in the orthogonal complement to the subspace of constants
 in $L^2_{\mu}(X)$ on which $V$ acts isometrically.

 \begin{proposition}{\it A Markov operator $V_{\Pi}$ is totally nonisometric if
 and only if $\Pi$ is a prime polymorphism.}
 \end{proposition}

 The dual notions of coisometrical and noncoisometrical Markov
 operators and connections with coprime and noncoprime
 polymorphisms are defined in a natural way.

 A mixing Markov operator is totally nonisometric, but we are interested
 in Markov operators that are far from isometries (in other words,
 in polymorphisms that are far from automorphisms)
 and far from mixing ones. Examples
 of prime nonmixing polymorphisms and, equivalently, totally nonisometric
 nonmixing operators play the key role in our theory; the existence
 of totally nonisometric nonmixing Markov operators is not  {\it a priori}
 obvious.

 In the terminology of the book \cite{FN}, a totally nonisometric
 nonmixing Markov operator is a Markov operator of type $C_{1,1}$
 (for the one-sided case, $C_{1, \cdot}$ or  $C_{\cdot,1}$).

 The first example of this type was given in \cite{R}; then in
 \cite{Vgrub} we suggested a general approach related to
 hyperbolic transformations, which we will use here (see also
 \cite{V1}). We will return to this in Sec.~3.

\section{Metric hyperbolic structure}

\subsection{Hyperbolic structure of a measure-preserving automorphism}
In this section, we formulate the main definitions.

 Suppose that $T$ is an ergodic measure-preserving automorphism of a Lebesgue space
 $(X,\mu)$. The following condition on the automorphism plays the key role
 in our considerations.

\smallskip\noindent
\textbf{Condition H.} There exists a $T$-invariant ergodic
equivalence relation $\chi$ and a polymorphism $\Pi$ that can be
represented as $$\Pi=\Phi \cdot T,$$ where  $\Phi$ is a nondegenerate
polymorphism associated with the
partition $\chi$;  at the same time, the following limits (in the weak
topology on the semigroup of polymorphisms) exist:
 \begin{equation}
\Lambda=\lim_{n\to \infty}\Pi^n\cdot T^{-n} \label{1}
\end{equation}
and
\begin{equation}
  \Gamma=\lim_{n\to \infty} T^{-n}\cdot\Pi^n.
 \label{2}
 \end{equation}
Besides, both polymorphisms $\Lambda$ and $\Gamma$ are
dense.

\begin{definition}
{\it A proper hyperbolic structure for an ergodic automorphism $T$ is
an ergodic equivalence relation $\chi$ for which Condition
{\bf H} (the existence of a polymorphism $\Pi$, etc.) holds.

An automorphism $T$ for which there exists at least one hyperbolic
structure will be called a hyperbolic automorphism; in this case,
we will say that $\chi$ is a homoclinic equivalence relation
\footnote{The notion of homoclinic equivalence relation for
automorphisms was introduced and used by M.~I.~Gordin \cite{Gord}
for other purposes. His definition is different, and we will
discuss its connections with our definition below and elsewhere.}
for the automorphism $T$. The same partition $\chi$ defines a
hyperbolic structure for the automorphism $T^{-1}$.}
\end{definition}

     The convergence of the infinite products in formulas (1) and (2) above is
 the main condition of our construction; in a sense, it is
 equivalent to the existence of stable and unstable
 foliations in the classical smooth theory of hyperbolic systems. It is
 easy to check that the polymorphisms $\Lambda$ and $\Gamma$ are also
 associated with the partition $\chi$.
 Note that for a given hyperbolic
 structure $\chi$,
 the choice of a polymorphism $\Phi$ and, consequently,
 of a polymorphism $\Pi$
 satisfying Condition {\bf H} is not unique;
 of course, both limits $\Lambda$ and $\Gamma$
 depend on the choice of $\Phi$.

 Nevertheless, technically,
 the central role is played by the polymorphism $\Phi$ associated with the
 relation $\chi$; therefore, we will rewrite the above limits
  in several forms, using the polymorphism $\Phi$
  instead of $\Pi$. Let
 $\Phi_k=T^k\Phi T^{-k}$, where $k\in \mathbb Z$, and $\Phi \equiv
 \Phi_0$; then we can rewrite these limits as
\begin{equation}
 \Lambda=\lim_{n\to \infty}\Phi\cdot
T\Phi T^{-1}\ldots T^n\Phi T^{-n},\label{3}
\end{equation}
or
$$
\Lambda = \lim_{n\to \infty} \prod_0^n \Phi_k = \prod_0^{\infty}\Phi_k; \eqno(3')
$$
analogously,
\begin{equation}
\Gamma=\lim_{n\to \infty}T^{-n}\Phi T^n \cdots T^{-1}\Phi T\cdot
\Phi,\label{4}
\end{equation}
or
$$\Gamma = \lim_{n\to \infty} \prod_{-n}^0
   \Phi_k = \prod_{-\infty}^0 \Phi_k. \eqno (4')
$$
Thus  we obtain the following proposition.

\begin{proposition}
{The convergence of two products in {\rm ($3'$)} and {\rm ($4'$)}
to dense polymorphisms $\Lambda$ and $\Gamma$
     (together with
the condition that $\Phi$ is a nondegenerate polymorphism associated
with the $T$-invariant  partition $\chi$) is equivalent to
Condition {\bf H}.}
\end{proposition}

In particular, we have the following important corollary.

 \begin{corollary} The sequence of polymorphisms $\Phi_k=T^k\cdot \Phi
 \cdot T^{-k}$ weakly tends to the identity automorphism as $|k| \to
 \infty$:
 $$\lim_{|k|\to \infty}T^k\cdot \Phi \cdot T^{-k}={\rm Id}. \eqno (\diamond)$$
\end{corollary}

The question is what rate of convergence can have the
 left-hand  side of $(\diamond)$ for various examples.

\subsection{Quasi-similarity of automorphisms and polymorphisms.}

The most essential ingredient of our construction is the
polymorphism $\Pi$.

\begin{theorem}
   Under Condition {\bf H}, the following formulas hold:
   $$\Pi \cdot\Lambda = \Lambda\cdot T \eqno (5)$$
   and
   $$\Gamma\cdot \Pi = T \cdot\Gamma. \eqno (6)$$
If relations {\rm (5)}, {\rm (6)} hold,
then we can claim that the automorphism $T$
 is quasi-similar to the polymorphism $\Pi$.
\end{theorem}

\begin{proof}
Because of the importance of equations (5),(6),
we present some calculations.
   Using our notation and the above formulas, we can rewrite these equations
   as
\begin{eqnarray*}
   \Pi\cdot \Lambda&=&  \Phi_0 \cdot T \cdot \lim_{n\to \infty} \prod_0^n \Phi_k
   =\lim_{n \to \infty}[\Phi  T \cdot \Phi \cdot
    T \Phi  T^{-1}\dots\ldots T^n\Phi T^{-n}]\\
&=&\lim_{n \to \infty}[\Phi \cdot T\Phi T^{-1} \cdot T^2 \Phi T^{-2}\dots T^{n+1}\Phi
   T^{-(n+1)}]\cdot T = \lim_{n\to \infty}\prod_{k=0}^{n+1} \Phi_k \cdot T\\
&=&\Lambda \cdot T;
\end{eqnarray*}
in a shorter form,
    $$\Pi \cdot \Lambda  =\Pi \cdot \lim_{n \to \infty} \Pi^n T^{-n}
   =\lim_{n \to \infty} [\Pi^{n+1}T^{-(n+1)}] T=\Lambda \cdot T.$$
Similarly,
   $$\Gamma \cdot \Pi =\lim_{n \to \infty} \prod_{-n}^0 \Phi_k \cdot \Phi_0 \cdot T =
   T \cdot \lim_{n \to \infty}\prod_{-(n+1)}^0 \Phi_k=T \cdot \Gamma,$$
   or
   $$\Gamma \cdot \Pi = \lim_{n \to \infty}[T^{-n} \cdot \Pi^n]\cdot \Pi = T \lim_{n \to \infty}
   T^{-(n+1)} \cdot \Pi^{(n+1)} = T \cdot \Gamma.$$
   The theorem follows from these equations and the definition of
quasi-similarity,
   together with the above conditions on the density of the
   polymorphisms $\Lambda$ and $\Gamma$.
\end{proof}

Now we can refine the properties of the polymorphism $\Pi$.

\begin{theorem}
A polymorphism $\Pi$ that satisfies Condition {\bf H} is prime,
nonmixing, and non-co-mixing. More exactly, if a polymorphism
$\Pi$ satisfies relations {\rm (1)} and {\rm (2)} with some
ergodic automorphism $T$ and dense polymorphisms $\Lambda$ and
$\Gamma$, then it is prime, coprime, nonmixing, and noncomixing.
\end{theorem}

\begin{proof}
 First we will prove that $\Pi$ (or $\Pi^*$) is prime.
 Suppose that $\Pi$ is not prime; this means that there exists a
 nontrivial $\Pi$-invariant {\it measurable} partition $\zeta$ of
 $(X,\mu)$. We have $T^{-n}\Pi^n \zeta=T^{-n}\zeta$, whence $\gamma
\equiv\Gamma \zeta=\lim_{n\to \infty}T^{-n}\zeta$; the existence of
 this limit is possible only if the partition $\zeta$ is $T$-invariant
 and, consequently, $\Gamma$-invariant. Therefore, we can consider the
 actions of the automorphism $T_{\zeta}$ and the endomorphism $\Pi_{\zeta}$
 on the quotient space $X/\zeta$. Thus we reduce the problem
 to the following one.

 \begin{proposition}
 {If the sequence of products $R^n\cdot S^{-n}$,
 where $R \ne {\rm Id}$ is an endomorphism and $S\ne {\rm Id}$ is an automorphism,
 tends to some limit in the weak topology, then  $R=S$.}
 \end{proposition}

 \begin{proof}
 {Consider the endomorphism $Q=S^{-1}R$;
 then (as above) we have
 $$
 R^n \cdot S^{-n}= (QS)^nS^{-n}=
 Q\cdot SQS^{-1}\cdot S^2QS^{-2}\ldots S^{n-1}QS^{-n+1}\cdot S^{-1}.
 $$
 The existence of the limit means that the following weak limit
 also exists and is equal to the identity:
 $$\lim_{n \to \infty} S^n\cdot Q\cdot S^{-n}={\rm Id}.$$
 But if $Q$ is an endomorphism, this can happen only if $S={\rm Id}$,  which is
 not the case, or if $Q={\rm Id}$, which means that $R=S$.}
\end{proof}

The claim of the proposition is not true if $R$ is a polymorphism.

 Now suppose that $\Pi$ is mixing. Recall (see \cite{V1, V2}) that
 each polymorphism naturally defines a Markov chain.
 We use the following observation (see \cite{V2}).

 \begin{proposition}
 The shift in the space of realizations of the Markov chain corresponding
 to a prime mixing polymorphism $\Pi$ is a $K$-automorphism, and
 the Markov generator is a $K$-generator.
\end{proposition}

Consequently, if $\Pi$ is mixing, then, in view of the $K$-property,
the above limit is again the zero polymorphism:
$$\lim_{n \to \infty}\Pi^n \cdot T^{-n}=\Theta;$$
but this limit is equal to $\Lambda$, which is impossible. Thus
$\Pi$ is not mixing; the same is true for $\Pi^*$.
\end{proof}

Note that the assertion converse to that of the proposition is
also true, so this gives a criterion of $K$-generators of Markov
chains. The question what prime nonmixing and non-co-mixing polymorphism
$\Pi$ defines a hyperbolic structure with given polymorphism $\Pi$
requires more information on the properties of the Markov process
generated by the polymorphism; we consider the corresponding
construction elsewhere (see also \cite{V1}).

 \subsection{Classical examples of hyperbolic structures}

\begin{theorem} A smooth hyperbolic transformation of a compact manifold
with finite invariant measure (Anosov system with discrete time)
has a natural proper hyperbolic structure in the above sense.
\end{theorem}

\begin{proof} {Let $T$ be an Anosov transformation of a smooth compact
manifold with an invariant measure; as an ergodic equivalence
relation $\chi$ from the definition above, we choose the ordinary
homoclinic equivalence relation: two points are equivalent
if they belong to the same stable and unstable leaves.
In the algebraic case --- that of a hyperbolic automorphism of the torus ${\mathbb
T}^n$ --- the homoclinic partition is the orbit partition of the
action of ${\mathbb Z}^{n-1}$, the Dirichlet group. It is an ergodic
relation, because the corresponding partition has a trivial measurable
hull. As a polymorphism $\Phi$, we can take a polymorphism for
which the conditional measure $\nu^x$ at a point $x$ is a
nondegenerate measure concentrated on a finite subset
of the set of points homoclinic to the point $x$ ---  such a polymorphism
$\Phi$ is associated
       with $\chi$ in the sense of our definition.
Thus the polymorphism $\Pi$
sends a point $x$ to the homoclinic class of the point $Tx$. In
order to prove that such a polymorphism $\Phi$ exists, or that one can
find a measurable map $x \to \nu ^x$, it suffices to
choose two different measurable maps on the manifold, each
associating with every point $x$ a point $y(x)\ne x$
homoclinic to $x$. The more serious part of the proof, the
existence of the limits (1) and (2), or the existence of
polymorphisms $\Lambda$ and $\Gamma$ above, was given in
\cite{Vgrub, V1}; see also Sec.~3.6.}
\end{proof}

\smallskip\noindent
{\bf Condition ($*$).} Let $T$ be a $K$-automorphism, and let $\xi$ be
a finite or countable $K$-generator of $T$ satisfying the following additional
property:
$$\bigwedge_{n=0}^{\infty}\bigvee_{|k|>n}
T^k\xi =\nu; \eqno (*)$$
here $\nu$ is the trivial measurable
partition.

 It is well known that not all $K$-generators of a
$K$-automorphism, and even of a Bernoulli automorphism, satisfy
property $(*)$.

\begin{theorem}
{If a $K$-automorphism satisfies condition {\rm ($*$)}, then it is
hyperbolic.}
\end{theorem}

\begin{proof}
 {Define an equivalence relation $\chi$ (= partition) as the
nonmeasurable partition obtained as the set-theoretic intersection
of the partitions from the previous expression:
$$\chi = \bigcap_{n=0}^{\infty}\bigvee_{|k|>n} T^k\xi. \eqno (**)$$
If we realize the automorphism $T$ as the right shift in the space of
two-sided sequences (states of the process), then two sequences
$x=\{x_i\}$, $y=\{y_i\}$, $i \in \mathbb Z$, belong to the same element of
$\chi$ if there exists $k \in \mathbb N$ such that $x_i=y_i$ for each $i>|k|$.
As the limit of a decreasing sequence of measurable
partitions with finite or countable blocks, $\chi$ is a hyperfinite
(or tame) partition. Condition $(*)$ means that the equivalence
relation $\chi$ is ergodic. A direct construction of
polymorphisms $\Phi$ and $\Pi$ with required properties is given in
Sec.~3.6. Thus $\chi$ determines a proper hyperbolic structure
for the automorphism $T$.}
\end{proof}

\smallskip\noindent
{\bf Remark.} For all $K$-automorphisms known at present (2005)
there exists a $K$-generator satisfying property $(*)$.\footnote{The author
is grateful to Professor J.-P.~Thouvenot for this information.} The
open question is whether such a generator exists for all
$K$-automorphisms.

\subsection{Geometrical interpretation}

 Relations (5) and (6) (together with (1) and (2))
 have an important {\it geometrical interpretation}. We interpret
 conditions (5) and (6); below the measure $\mu_{\Psi}^z$ is the
$\Psi$-image of a point $z$ (see definitions):
$$
\mu_{\Lambda}^{Tx}(\cdot) = \int \mu_{\Pi}^y(\cdot)d\mu_{\Lambda}^x(y),
\eqno(5')
$$
$$
\mu_{\Gamma}^y(T\cdot)= \int\mu_{\Lambda}^x(\cdot)d\mu_{\Gamma}^y(x).\eqno (6')
$$
But the action of a polymorphism on a given space
can be naturally extended to the action on probability measures on the same
space (convolution); using this, we can rewrite the formula as
follows:
 $$\mu_{\Lambda}^{Tx}(\cdot)=(\Pi*\mu_{\Lambda}^x)(\cdot)$$
and, respectively,
 $$\mu_{\Gamma}^y(T\cdot)=(\mu^y_{\Gamma}*\Pi)(\cdot).$$
These equalities show that the right (resp., left) action of the
automorphism $T$ on the set of measures $\mu_{\Lambda}^x$
(resp., $\mu_{\Gamma}^x$), $x \in X$,
is the same as the left (resp., right)
action of the polymorphism $\Pi$ on these sets of measures. This
is an explanation of the quasi-similarity between the automorphism $T$
and the polymorphism $\Pi$.

The polymorphisms $\Lambda$ and $\Gamma$ are of special
interest from the point of view of the corresponding Markov processes;
see \cite{V1} and \cite{Vgrub}.

An ergodic automorphism can have several metrically
nonisomorphic hyperbolic structures or have no such structures.
The main problem is to characterize automorphisms that have
hyperbolic structures and to classify these structures;
this problem is new and
has no answer up to now. We will present some results in this
direction.

\bigskip

The definitions of hyperbolic structures,
quasi-similarity, and other notions discussed above
could be easily
reformulated in terms of unitary and Markov positive operators
in the space $L^2$ (see \cite{V2}). We restrict ourselves only to an
operator reformulation of quasi-similarity.

Let $U_T$ be the unitary operator in $L^2$ corresponding to a
measure-pre\-serv\-ing automorphism $T$, and let $V_{\Pi}$ be the Markov operator
corresponding to a polymorphism $\Pi$ (see definitions). Then, under
conditions (3) and (4), the following limits of Markov operators
in the weak operator topology do exist:
\begin{eqnarray*}
\lim_{n \to \infty} V_{\Pi}^n U_T^{-n}&=&L,\\
\lim_{n \to \infty}U_T^{-n}V_{\Pi}^n&=&G,
\end{eqnarray*}
and the operators $U_T$ and $V_{\Pi}$ are quasi-similar:
$$V_{\Pi}L=LU_T, \qquad GV_{\Pi}=U_TG.$$

Recall that the polymorphism $\Pi$ is nonmixing (non-co-mixing) and prime
and, consequently, the Markov operator $V_{\Pi}$ is totally nonisometric and
nonmixing (non-co-mixing), so this is an example of
quasi-similarity between a unitary and a totally nonisometric operator.
This is a positive analog of operators of class $C_{1,1}$ in the
sense of \cite{FN}. The existence of such examples is not obvious.

\subsection{Left and right semi-hyperbolic structures}

Now we define  structures that are weaker than the hyperbolic one.
Let $T$ be an ergodic measure-preserving automorphism. We define
left and right semi-hyperbolic structures.

\smallskip\noindent
\textbf{Condition SH.} There exists a $T$-invariant ergodic
equivalence relation $\chi_l$ (resp., $\chi_r$) and a
polymorphism $\Pi_l$ (resp., $\Pi_r$) that can be represented
as $$\Pi_l=\Phi_l \cdot T$$ (resp., $$\Pi_r=\Phi_r\cdot T),$$
where the polymorphism $\Phi_l$ (resp., $\Phi_r$) is associated with
the partition $\chi_l$ (resp., $\chi_r$) and such that the following
limits in the weak topology on the semigroup of polymorphisms exist:
$$ \Lambda_l=\lim_{n\to \infty}\Pi_l^n\cdot T^{-n} \eqno(1_l)$$
(resp.,
$$\Gamma_r=\lim_{n\to \infty} T^{-n}\cdot\Pi_r^n.) \eqno(2_r)$$
Besides, both polymorphisms $\Lambda_l$ and $\Gamma_r$ are
dense.

As in Definition 4, the polymorphism $\Pi_l$ (resp., $\Pi_r$)
is prime and nonmixing (resp., non-co-mixing).

\begin{definition}
{An ergodic equivalence relation $\chi_l$ (resp., $\chi_r$) defines
a left (resp.,  right) semi-hyperbolic structure of an automorphism
$T$ if condition $(1_l)$ (resp., $2_r$) of Condition {\bf SH}
holds. If both conditions hold for some polymorphisms $\Pi_l$
and $\Pi_r$, we will say that they define a semi-hyperbolic
structure for the automorphism $T$.}
\end{definition}

Similarly to Theorem 1, we have the following result.

\begin{theorem}
The following relations hold:
\begin{eqnarray*}
\Pi_l\cdot\Lambda_l&=&\Lambda_l\cdot T,\\
\Gamma_r \cdot\Pi_r&=&T\cdot\Gamma_r.
\end{eqnarray*}
\end{theorem}

Consequently, the automorphism $T$ is a quasi-image of the polymorphism
$\Pi_l$, and the polymorphism $\Pi_r$ is a quasi-image of the automorphism
$T$.

\medskip
 If $T$ is the shift in the space of realizations of a stationary
 random process, one can take as $\chi_l$ or $\chi_r$ the
 partitions with fixed past or future of this process.
 In smooth hyperbolic theory, the partitions $\chi_l$ and
 $\chi_r$ can be chosen to be the stable and unstable foliations, respectively;
 in that theory, they exist simultaneously.  In
 general, in the above definition there are no connections
 between the left and right semi-hyperbolic structures. Thus we can
 consider many variants and examples of semi-hyperbolic structures
 in the sense of our definition. If the left and right structures agree, in
 the sense that the supremum of two relations $\chi_l$ and $\chi_r$ (or,
 in terms of partitions,
 the product of the partitions $\chi_l$ and $\chi_r$)
 is an ergodic relation (resp., an ergodic partition),
 then we have a {\it proper hyperbolic structure} in the sense of the
 main definition of Sec.~3.1, and this product is a homoclinic partition.

\subsection{Semi-hyperbolic structure of $K$-automorphisms}

  We will prove that every $K$-automorphism possesses a
  semi-hyperbolic structure.

\begin{theorem}
   {For every $K$-automorphism $T$ there exists a semi-hyperbolic
   structure, i.e., there exists a prime nonmixing polymorphism
   $\Pi_l$ that defines a left semi-hyperbolic structure of the automorphism
   $T$ and a prime non-co-mixing polymorphism $\Pi_r$ that defines a
   right semi-hyperbolic structure.}
\end{theorem}

\begin{proof}
      It is clear from the definition that if $T$ has a left (right)
  semi-hyperbolic structure, than it is a right (left) semi-hyperbolic
  structure for the automorphism $T^{-1}$. Now if $T$ is a
$K$-automorphism,
  then $T^{-1}$ is also a $K$-automorphism. Thus it suffices
  to prove that every $K$-automorphism has a left semi-hyperbolic structure.

We will construct series of polymorphisms that will
be random perturbations of special type of the initial
$K$-automorphism.
Our construction is a detailed version of the previous examples from
the papers  \cite{Ros, Vgrub, V1}.

Let $T$ be an arbitrary $K$-automorphism. By well-known theorems
(see, e.g., \cite{KSinF}), we can realize $T$
 as the right shift in the space $\cal X$ of two-sided sequences in a finite
 or countable alphabet $A$; the space $\cal X$ is equipped with
 a shift-invariant measure $\mu$ and has trivial (in the sense
 of the measure $\mu$) tail algebras in the past and in the future.
 This means that if we consider the one-sided right shift in the space
 of one-sided sequences ${\cal X}=\{\{x_n\}_{n<0}\}$ (we denote it by
 the same letter), then we have a decreasing sequence of measurable
 partitions $\zeta_n\equiv T^{-n}\varepsilon $, and for every $n$
 the partition $\zeta_n$ has countable or finite fibers; here $T$
 is the right shift in the space ${\cal X}$ and $\varepsilon$ is
 the partition  of the space $\cal X$ into separate points.
  We have $\bigwedge_n \zeta_n =\nu$, where $\nu$ is the trivial partition.
  Our first goal, according to the definition of a left hyperbolic structure,
  is to construct a measure-preserving polymorphism
  $\Phi_l=\Phi$ that acts in the space $\cal X$ and has the
  following structure: for almost all $x$, the supports of the
  measures $\Phi(x)\equiv\mu_x$ belong to the element $\zeta_n (x)$
  of the partition $\zeta_n$ for sufficiently large $n$.
    We have a space $\cal X$ with measure $\mu$ and a decreasing sequence of
    measurable partitions $\zeta_n$ that tends to the trivial partition. The
    desired polymorphism $\Phi$ (and later $\Pi$) will be a kind of random
    walk over several automorphisms with quasi-invariant measure.

    Note that since the intersection $\bigwedge_n \zeta_n$ is trivial,
    the number of points in the elements of $\zeta_n$ tends to infinity
    (or already equal to infinity). Thus
    for every small $\delta>0$ there exist a positive integer $n_{\delta}\in \mathbb
    N$ and a measurable set $C_{\delta} \subset X$ of $\mu$-measure
    greater than $1-\delta$ such that for every point $x \in C_{\delta}$,
    the element of the partition $\zeta_n$ containing $x$ contains at least
    {\it four} different points (including $x$).
    Let us take a measurable refinement $\eta_n$ of the restriction of $\zeta_n$
    to the set $C_{\delta}$ with all elements consisting of four
    points. If the number of points in the elements is not divisible
    by four, we form the set $C'_{\delta}$ from all the remaining points and
    join it to the set $\cal X \smallsetminus C_{\delta}$.
    Then we restrict our partitions $\zeta_n$ with $n >n_{\delta}$ to the set
    $(\cal X \smallsetminus C_{\delta})\cup C'_{\delta}$
    and repeat this procedure again. Finally, we obtain a measurable
    partition $\eta$ of the whole space $\cal X$, and
    each of its elements is a {\it refinement} of some element of the
    partition $\zeta_n$ for some $n$.

    Each element of the partition $\eta$ consists of four points, and we label
    them with the numbers $1$, $2$, $3$, $4$ in a
    measurable way (so that the set of points with label $i$ ($i=1,2,3,4$)
    is measurable). For
    every point $y$, denote by $n(y)$ the label of $y$.
    Define three involutions on four points $v^1,v^2,v^3$
    as follows:
    $v^1=(1,2)(3,4)$, $v^2=(1,3)(2,4)$, $v^3=(1,4)(2,3)$.
    Denote by $p_i=p_i(y)$ the conditional measures of points in the element
    of the partition $\eta$ that contains $y$.
    Using the combinatorial lemma given below, we introduce
    a polymorphism $\Phi$ as follows:
    $$\Phi(y)=v^i(y) \qquad\mbox{with probability} \quad q_{i,n(y)}.$$

    Note that the partition $\eta$ is, by definition, a fixed partition
    for the polymorphism $\Phi$. Thus the
    images of $y$ under $\Phi$ are points  (not equal to $y$) of the same element
  of $\eta$ that contains $y$ with
  some probabilities.
  The fact that $\Phi$ preserves the measure follows from the
  construction of $q_{i,n(y)}$; we must only
  mention that because of the transitivity of conditional measures,
  the conditional measure of the point $x^i$ with respect
  to the partition $\eta$ is the same as
  the conditional measure with respect to the
  partition $\zeta_n$.

  Now we define a polymorphism $\Pi$ on the space
  $\cal X$ as follows: $\Pi=\Phi \cdot T$.

  Since $\Pi$ is the product of two measure-preserving
  automorphisms $\Phi$ and $T$,
  it also preserves the measure $\mu$.

     From the definition of $\Pi$  we see that
     it sends a sequence $x=\{x_n\}$
     to the shifted sequence $Tx$ and then changes at
     random a {\it finite  number of digits}.
     The fact that the polymorphism changes
     a finite number of coordinates follows
     from the fact that the elements of the partition $\eta$
     (in which the involutions $v^i$ act) are contained
     in some element of the partition
     $\zeta_n$ for some $n$, so at most $n$ coordinates can
     be changed. Of course, this $n$ depends on $x$, thus it can be
     arbitrarily large.

      In order to finish the proof, we need to prove that

     1) the polymorphism $\Pi$ is prime;

     2) there exists $\lim_{n\to \infty}\Pi^n \cdot T^{-n}$.

     The primality (the absence of nontrivial invariant measurable partitions)
           follows from the fact that, by definition, an invariant partition
     for $\Pi$ that does not coincide with $\varepsilon$
     must be less (coarser)
     (for definitions, see, e.g., \cite{V3})
     than the intersection $\bigwedge_n \zeta_n$,
     which is the trivial partition $\nu$.

     The existence of the limit follows from the structure of
     $\Pi$. Indeed, the polymorphism $\Pi^n$
     shifts every sequence by $n$ and changes at random a
     finite number of digits so that no
     digit is changed infinitely many times.
     Thus  $\Pi^n \cdot T^n$ is a polymorphism that for every
     $x$ changes finitely many digits of $x$, and the coordinates
     of these digits go to infinity, so each coordinate stabilizes,
     and the corresponding measures converge.
\end{proof}

    Now we formulate a simple combinatorial lemma that we have used
    in the proof of the theorem.

\begin{lemma}
Let $p_i$, $i=1,2,3,4$, be an arbitrary probability vector of
length four. There exists a matrix $\{q_{i,j}\}_{i,j=1}^4$ with
$q_{i,i}=0$, $q_{i,j}\geq 0$, $i,j=1,2,3,4$, with given
marginal projections: $\sum_i q_{i,j}=p_j$, $\sum_j q_{i,j}=p_i$.
\end{lemma}

   The proof of this lemma is straightforward.
   We may say that the matrix $\{q_{i,j}\}$ determines a measure-preserving
   polymorphism of the space.
   Since the lemma is valid for any number of
   points greater than four, the partition $\eta$
     can also be chosen with arbitrarily many points.

     For $K$-automorphisms satisfying property ($*$) (see Sec.~3.3), we can
     choose a polymorphism $\Pi$ that  simultaneously defines a
     left and right semi-hyperbolic structures and, consequently,
     a hyperbolic structure. The supports of the conditional measures
     are the blocks of the homoclinic partition. As we have already mentioned,
     all known $K$-automorphisms have this partition.
     But for a general $K$-automorphism the situation is unclear.

\subsection{A conjecture and a problem}
\smallskip\noindent
     {\bf 1. Conjecture.} Each hyperbolic automorphism
is a $K$-automorphism satisfying property ($*$); each automorphism that is
quasi-similar to a prime nonmixing and non-co-mixing polymorphism has a
hyperbolic structure defined by this polymorphism.

\smallskip\noindent
 {\bf  2. Problem.} Is property ($*$) equivalent to the
$K$-property? Or, is every $K$-automorphism hyperbolic?

\end{document}